\documentclass[12pt]{article}
\usepackage{amsmath,amssymb,amsfonts,color}
\usepackage[latin1]{inputenc}
\def\eop{\hfill\rule{2.5mm}{2.5mm}}
\def\pf{\par\smallbreak\noindent {\bf Proof.} \ }
 \newtheorem{thm}{Theorem}[section]
 
 \newtheorem{lem}[thm]{Lemma}

 \numberwithin{equation}{section}

\begin{document}

%
%
%
%
%
%
%
%
%

\begin{center}{\LARGE{\bf Traceability of positive integral operators in the absence of a metric}}\end{center}

\author[{Mario H. Castro, Valdir A. Menegatto and Ana P. Peron}



\footnote{First author partially supported by CAPES}\footnote{ Second author partially supported by FAPESP $\#$ 2010/19734-6}


\date{January 1, 2004}

\begin{abstract}
We investigate the
traceability of positive integral operators on $L^2(X,\mu)$ when $X$
is a Hausdorff locally compact second countable space and $\mu$ is a
non-degenerate, $\sigma$-finite and locally finite Borel measure.\
This setting includes other cases proved in the literature, for
instance the one in which $X$ is a compact metric space and $\mu$ is
a special finite measure.\ The results apply to spheres, tori and
other relevant subsets of the usual space $\mathbb{R}^m$.
\end{abstract}


\section{Introduction}

Let $X$ be a Hausdorff locally compact and second countable
topological space endowed with a non-degenerate, $\sigma$-finite and
locally finite Borel measure $\mu$.\ In this paper, we shall
investigate the traceability of integral operators $\mathcal{K}:L^2(X,\mu)
\rightarrow L^2(X,\mu)$ generated by a suitable kernel $K:X\times
X\to \mathbb{C}$ from $L^2(X\times X,\mu\times\mu)$.\ The title of the paper
refers to the fact that the space $X$ carries a topological
structure rather than a metric one.\  The setting just described
allows the space $L^2(X,\mu)$ to have a countable complete
orthonormal subset ([7, p.92]) while the operator $\mathcal{K}$, which is
given by the formula
\begin{equation}\label{1.1}
\mathcal{K}(f):=\int_X K(\cdot,y)f(y)\, d\mu(y),\quad f\in L^2(X,\mu),
\end{equation}
becomes compact.\ As so, the spectral theorem for compact operators
is applicable and $\mathcal{K}$ can be represented in the form
\begin{equation}\label{1.2}
\mathcal{K}(f)=\sum_{n=1}^{\infty}\lambda_n \langle f,f_n\rangle f_n, \quad f \in L^2(X,\mu),
\end{equation}
in which $\{\lambda_n\}$ is a sequence of real numbers (possibly finite)
converging to 0 and $\{f_n\}$ is a complete orthonormal sequence in
$L^2(X,\mu)$.\ The symbol $\langle\cdot,\cdot \rangle$ will stand for the
usual inner product of $L^2(X,\mu)$.

The basic requirement on the kernel $K$ will be its positive
definiteness.\ A kernel $K$ from $L^2(X\times X,\mu\times\mu)$ is
$L^2(X,\mu)$-{\em positive definite} when the corresponding integral
operator $\mathcal{K}$, is positive:
\begin{equation}\label{1.3}
\langle\mathcal{K}(f),f\rangle \geq 0, \quad f \in
L^2(X,\mu).
\end{equation}
Fubini's theorem is all that is need in
order to show that a $L^2(X,\mu)$-positive definite kernel is
hermitian $\mu\times\mu$-a.e..\ As so, the integral operator $\mathcal{K}$ is
automatically self-adjoint with respect to $\langle\cdot,\cdot \rangle$.\
In particular, the sequence $\{\lambda_n\}$ mentioned in the previous
paragraph needs to be entirely composed of nonnegative numbers.\ In
the present paper, we shall assume they are listed in a decreasing
order, with repetitions to account for multiplicities.

Under the conditions established above, the specific aim of this
paper is to establish additional conditions on $K$ in order that
$\mathcal{K}$ be {\em trace-class}, that is,
\begin{equation}\label{1.4}
\sum_{f\in\mathfrak{B}}\langle\mathcal{K}^*\mathcal{K}(f),f\rangle^{1/2}<\infty
\end{equation}
for every orthonormal basis $\mathfrak{B}$ of $L^2(X,\mu)$.\ In the
formula above, $\mathcal{K}^*$ is the adjoint of $\mathcal{K}$.\ We refer the reader
to $[5,10,11]$ for more information on trace-class operators.

The main result in this paper can be seen as a generalization of
another one originally proved in [11] for the case $X=[a,b]$.\ The
proof there used in a key manner the so-called Steklov's smoothing
operator to construct an averaging process to generate a convenient
approximation to $\mathcal{K}$.\ The upgrade to the case in which $X$ is a
subspace of $\mathbb{R}^n$ was discussed in [8] and references therein.\ By
assuming that the Lebesgue measure of nonempty intersections of $X$
with open balls of $\mathbb{R}^n$ was positive and using auxiliary
approximation integral operators generated by an averaging process
constructed via the Hardy-Littlewood theory, the main result in [8]
described necessary and sufficient conditions for the traceability
of the integral operator, under the assumption of positive
definiteness of the kernel.\ The process used in [8] and other references as well
provides a way to deal with the generating kernel on the diagonal of $X\times X$
and it is convenient when the kernel is not continuous.\ Despite using a similar average process,
another achievement in the present paper is the inclusion of a setting in which the measure does not need to be finite.

Since our spaces are no longer metric, the Hardy-Littlewood theory
in the average arguments need to be replaced or adapted.\ We will
use techniques involving the construction of auxiliary integral operators based on martingales
constructed from special partitions of $X$, following very closely the
development of Brislawn in [2].\ A similar construction have appeared in
[6] in an attempt to generalize Brislawn results to $L^p$ spaces.\ The main difference between the
construction to be delineated here and those in [2] and [6] is that, in the present
one, we need to guarantee that the elements in the partitions belong
to the topology of $X$.\ This is the exact point where the
assumption of local compactness will play an important role.

For the sake of completeness we mention references $[1,14]$ where
other characterizations for traceability were obtained.

An outline of the paper is as follows.\ Section 2 contains the basic
information on martingales used in the paper, along with the key
construction we will need in order to introduce approximating
auxiliary operators in Section 3.\ There, the main technical results
are established and proved.\ Section 4 contains the main results of
the paper, including a convenient equivalence for traceability.

\section{A special martingale}

This section contains several results involving a special martingale
on $X$.\ Some of them are just refined versions of results described
in Section 2 of [2].\ However, the reader is advised that the basic
references we used for the concepts and results either quoted or
used here are $[4,15]$.

Let $(X,\mathcal{M},\sigma)$ denote a $\sigma$-finite measure space and $\mathcal{F}$
a sub-$\sigma$-algebra of $\mathcal{M}$ for which $(X,\mathcal{F},\sigma)$ is a
$\sigma$-finite measure space too.\ If $f:X \rightarrow \mathbb{C}$
is $\mathcal{M}$-measurable, Radon-Nikodyn's theorem asserts that we can find
a unique $\mathcal{F}$-measurable function $g: X\rightarrow \mathbb{C}$ so
that
\begin{equation}\label{1.5}
\int_A f \, d\sigma= \int_A g \, d\sigma,\quad A\in\mathcal{F}.
\end{equation}
The function $g$ is called the {\em conditional expectation} of $f$
relative to $\mathcal{F}$ and is written $g=E(f|\mathcal{F})$.\ If $\{\mathcal{F}_n\}$ is a
family of sub-$\sigma$-algebras of $\mathcal{M}$, a sequence $\{f_n\}$ of
$\mathcal{M}$-measurable functions on $X$ is a {\em martingale} if every
$f_n$ is $\mathcal{F}_n$-measurable and $E(f_n|\mathcal{F}_m)=f_m$, $m<n$.

Next, we remind the reader about the basic setting we are assuming
in the paper: $X$ is a Hausdorff, locally compact and second
countable topological space endowed with a non-degenerate, locally
finite and $\sigma$-finite Borel measure $\mu$.\ In addition to
that, we will write $\mathcal{B}_X$ to denote the Borel
$\sigma$-algebra of $X$.

Invoking the first countability axiom, we may infer that every point
of $X$ possesses an open neighborhood.\ Since $X$ is Hausdorff and
locally compact, these neighborhoods can be assumed to be the
interior of a compact set.\ Thus, due to the local finiteness of
$(X,\mu)$, we can assume, in addition, that the
open neighborhoods of elements of $X$ have finite measure.

We intend to construct a special sequence of partitions of $X$ from
an open covering $\{\mathcal{A}_x\}_{x\in X}$ of it, composed of neighborhoods
of the type just described, and use them to define a particular
martingale.\ If such a covering has been fixed, Lindel\"{o}ff's
theorem ([13, p.191]) implies that we can extract from it a
countable sub-collection $\{\mathcal{A}_n\}$, still covering $X$.\ Such
sub-collection can be used in the construction of a first stage
partition $\mathcal{P}_0$ of $X$, following these steps: the first two
elements in the partition are $\mathcal{A}_0$ and its frontier
$\partial\mathcal{A}_0$.\ Observing that $\{\mathcal{A}_n\backslash\overline{\mathcal{A}_0}\}$ is an
open and countable covering of $X\backslash\overline{\mathcal{A}_0}$, we pick
$\mathcal{A}_1\backslash\overline{\mathcal{A}_0}$ and $\partial \mathcal{A}_1\backslash\overline{\mathcal{A}_0}$
to include in the partition.\ The family $\{\mathcal{A}_n\backslash
\overline{\mathcal{A}_0\cup\mathcal{A}_1}\}$ is an open and countable covering of
$X\backslash\overline{\mathcal{A}_0\cup\mathcal{A}_1}$.\ We proceed, including its elements
$\mathcal{A}_2\backslash\overline{\mathcal{A}_0\cup\mathcal{A}_1}$ and
$\partial\mathcal{A}_2\backslash\overline{\mathcal{A}_0\cup\mathcal{A}_1}$ in the partition.\
Proceeding inductively, we complete the construction of $\mathcal{P}_0$,
which is countable and entirely composed of Borel sets of finite
measure.\ Since Theorem 7.8 in [9] implies that $\mu$ is regular,
all the sets of the form
$\partial\mathcal{A}_n\backslash\overline{\mathcal{A}_0\cup\dots\cup\mathcal{A}_{n-1}}$ in $\mathcal{P}_0$ have measure zero.

In the next step, we construct a sequence $\{\mathcal{P}_n\}$ of
partitions of $X$ from $\mathcal{P}_0$, using as we can, a countable
basis $\{\mathcal{U}_n\}$ for the topology of $X$.\ For $n=0,1,\ldots$, we put
$$\mathcal{P}_{n+1}=\{\mathcal{U}_n\cap\mathcal{A}:\, \mathcal{A}\in\mathcal{P}_n\}\cup \{(X\setminus\overline{\mathcal{U}_n})\cap\mathcal{A}:\,
\mathcal{A}\in\mathcal{P}_n\} \cup \{\partial\mathcal{U}_n\cap\mathcal{A}:\, \mathcal{A}\in\mathcal{P}_n\}.$$
Clearly, $\mathcal{P}_{n+1}$ refines $\mathcal{P}_n$ and the sequence
$\{\mathcal{F}_n\}$ of the corresponding $\sigma$-algebras generated by those
partitions increases to $\mathcal{B}_X$.\ In addition, every $(X,\mathcal{F}_n,\mu)$
is $\sigma$-finite.

It is easy to see that for each $x\in X$ and each positive $n$,
there exists a unique set $O_n(x)\in\mathcal{P}_n$ such that $x\in O_n(x)$.\
We denote by $\mathfrak{N}$ the subset  of $X$ containing all $x\in
X$ for which $\mu(O_m(x))=0$, for some $m\geq 0$.\ Since the
sequence $\{O_{n}(x)\}$ is telescoping, the equality $\mu(O_m(x))=0$
implies $\mu(O_n(x))=0$, $n\geq m$.\ Being each $\mathcal{P}_n$ countable, it
is easily seen that $\mu(\mathfrak{N})=0$.

The very same arguments used in [15, p.89] show that for every $x\in
X\setminus\mathfrak{N}$ and every positive $n$,
the conditional expectation
$E_n(f)$ of $f$ relative to $\mathcal{F}_n$ is given by the formula
\begin{equation}\label{2.2}
E_n(f)(x)={1 \over \mu(O_n(x))} \int_{O_n(x)} f \,
d\mu.
\end{equation}
The sequence $\{E_n(f)\}$ defines a
martingale generated by just one (measurable) function, the {\em
martingale associated with $f$}.\ Examples related to constructions
similar to the one above can be found in [15, p.88].

The section will be completed with a list of results involving the
previous formula and the {\em maximal function} $Mf$ of the
martingale associated with $f$.\ Such martingale is defined by the formula

\begin{equation}\label{2.3}
Mf(x):=\sup\{|E_n(f)(x)|:n=1,2,\ldots\},\quad x\in X.
\end{equation}
Since the results are quite general and are not attached to the
particular setting introduced above, we will include sketches of the
proofs for the convenience of the reader.

A classical result concerning the maximal function ([15, p.91])
implies that if $p \in (0,\infty)$ then
\begin{equation}\label{2.4}
\| Mf \|_p\leq c_p \|f\|_p,\quad f\in L^p(X,\mu),
\end{equation}
 where
$c_p$ is a constant depending on $p$ only and $\|\cdot\|_p$ denotes
the usual norm of $L^p(X,\mu)$.\ As for the conditional expectation,
it transforms convergence in the mean into convergence $\mu$-a.e.\
Another basic result ([4, p.53] and [2, p.232]), commonly called
Doob's martingale convergence theorem, states that $E_n(f)$
converges to $f$ $\mu$-a.e., as long as $f\in L^p(X,\mu)$ and $p \in
[1,\infty]$.


Moving forward, the inequalities
\begin{equation}\label{2.5} |E_n(f)(x)|
\leq
Mf(x), \quad x \in X\setminus \mathfrak{N},\quad n\geq 1,
\end{equation} and
\begin{equation}\label{2.6}|f(x)|  \leq
       |f(x)-E_n(f)(x)|+
Mf(x),\quad x \in X\setminus \mathfrak{N},\quad n\geq 1,\end{equation}
are easily deducted.\ Combining the last one with Doob's martingale
convergence theorem, we are led to the inequality $|f|\leq Mf$,
$\mu$-a.e..\ As for the conditional expectation, we have the
following result found in [15, p.90]: if $p \in [1,\infty]$ and
$f\in L^p(X,\mu)$ then $\|E_n(f)\|_p\leq\|f\|_p$, $n\geq 1$.\ As a
consequence, the following theorem holds.

\begin{thm}
If $p \in [1,\infty]$ then the linear map
$E_n:L^p(X,\mu)\to L^p(X,\mu)$ is bounded.\ If $p=2$, then the
previous map is a self-adjoint operator.
\end{thm}

We close the section with a result for convergence in the mean of
the conditional expectation.

\begin{thm}If $f\in L^2(X,\mu)$ then $E_nf$ converges to
$f$ in the mean.\end{thm}

\pf If $g_n:=|f-E_n(f)|^2$, $n\geq 1$, the previous theorem yields
that $\{g_n\}\subset L^1(X,\mu)$.\ Now, inequality (2.5) leads to
\begin{equation}\label{2.7} |g_n(x)|\leq 2(|f(x)|^2+|E_n(f)(x)|^2)\leq 4 |Mf(x) |^2,\quad
x\in X\setminus \mathfrak{N},\quad n\geq 1.\end{equation} Clearly,
$Mf\in L^2(X,\mu)$ while Doob's convergence theorem gives us $g_n\to
0$ $\mu$-a.e..\ The dominated convergence theorem connects the final
arguments.\eop

\section{Approximating kernels}

This section is entirely composed of technical results involving a
family of operators constructed from the martingale defined in
Section 2.

Under the notation in Section 2, Theorem 7.20 in [9] informs that
the product measure $\mu\times\mu$ is a regular Borel measure on $X\times X$
and the sequence $\{\mathcal{P}_n\times\mathcal{P}_n\}$ of partitions of $X\times X$
increases to the Borel $\sigma$-algebra $\mathcal{B}_{X\times X}$ of
$(X\times X, \mu\times \mu)$.\ In particular, if $K\in
L^1_{loc}(X\times X,\mu\times\mu)$, the conditional expectation with
respect to the $\sigma$-algebra generated by the partition
$\mathcal{P}_n\times\mathcal{P}_n$ of $X\times X$ can be defined by the formula
\begin{equation}\label{3.1}
E_n(K)(u,v):={1 \over \sigma(O_n(u))\sigma(O_n(v))} \int_{O_n(u)}
\int_{O_n(v)} K(x,y)\, d\mu(y)d\mu(x).\end{equation}

Lemma 3.1 below provides information about a limit property
regarding the open sets $O_n(x)$ previously defined.\ We will use
the symbol $\chi_A$ to denote the characteristic function of the
subset $A$ of $X$.\ We remind the reader that given $x \in X$ and
$n\geq 1$, the construction introduced in the previous section shows
that there exists a unique $O_n(x) \subset P_n$ so that $x \in
O_n(x)$.

\begin{lem}If $x \in X$ and $n\geq 1$ then
\begin{equation}\label{3.2}\lim_{u\to
u_0}\chi_{O_n}(u)(x)=\chi_{O_n(u_0)}(x), \quad u_0\in X\setminus
\mathfrak{N}.\end{equation}\end{lem}

\pf Fix $x \in X$ and $n\geq 1$.\ If $u \in X$ then $x\in O_n(u)$ if
and only if $u\in O_n(x)$.\ Since
$\chi_{O_n(x)}(u)=\chi_{O_n(u)}(x)$, we can write
\begin{equation}\label{3.3}
|\chi_{O_n(u)}(x)-\chi_{O_n(u_0)}(x)|=|\chi_{O_n(x)}(u)-\chi_{O_n(x)}(u_0)|,\quad
u_0\in X. \end{equation} Next, if $u_0\in X\setminus \mathfrak{N}$,
the fact that $O_n(u_0)$ is open, leaves us with two cases: if $x\in
O_n(u_0)$ then $u_0\in O_n(x)$ and, at the limit, we can assume
$u\in O_n(u_0)=O_n(x)$ so that
\begin{equation}\label{3.4}
\lim_{u \to
u_0}|\chi_{O_n(x)}(u)-\chi_{O_n(x)}(u_0)|=|1-1|=0.\end{equation}
If $x\not\in O_n(u_0)$ then $u_0\not\in O_n(x)$, and assuming $u \in
O_n(u_0)$ as we can, we conclude that
\begin{equation}\label{3.5}
\lim_{u\to
u_0}|\chi_{O_n(x)}(u)-\chi_{O_n(x)}(u_0)|=|0-0|=0.\end{equation}
The proof is complete.\eop

It is now reasonable that the following result holds.

\begin{lem}If $u_0\in X\setminus \mathfrak{N}$ then $\lim_{u\to u_0}\mu(O_n(u))=\mu(O_n(u_0))$, $n=1,2,\ldots.$\end{lem}

\pf Since \begin{equation}\label{3.6}\mu(O_n(u))=\int_X \chi_{O_n(u)}(x)\, d\mu(x),\quad u\in
X,\end{equation} it follows that
\begin{equation}\label{3.7}
|\mu(O_n(u))-\mu(O_n(u_0))|\leq \int_X
|\chi_{O_n(u)}(x)-\chi_{O_n(u_0)}(x)|\, d\mu(x),\quad u \in
X.\end{equation}
As so, the assertion of the lemma will be proved if we can show that
\begin{equation}\label{3.8}
\lim_{u\to u_0}\int_X |\chi_{O_n(u)}(x)-\chi_{O_n(u_0)}(x)|\,
d\mu(x)=0, \quad u_0\in X\setminus \mathfrak{N}.\end{equation}
Hence, in view of the previous lemma, it suffices to show that the
integral and the limit in the previous equation commute.\ The family
$\{g_u\}$ defined by
\begin{equation}\label{3.9}g_u(x)=|\chi_{O_n(u)}(x)-\chi_{O_n(u_0)}(x)|, \quad u,x\in X,\end{equation} and the function
$g=\chi_{O_n(u_0)}$ belong to $L^1(X,\mu)$.\ Since $|g_u|\leq g$,
$\mu$-a.e., when $u \to u_0$, the desired commuting property follows
from the dominated convergence theorem.\eop

We now turn to kernels of the form
\begin{equation}\label{3.10}
D_n(u,x)= \frac{1}{\mu(O_n(u))}\chi_{O_n(u)}(x),\quad u,x\in X,\quad
n=1,2,\ldots.\end{equation}
and the corresponding integral operators $\mathcal{D}_n$ generated by $D_n$.\
For use ahead, we mention the immediate formula \begin{equation}\label{3.11}
E_n(\chi_{O_n(u)}\, f)=\mathcal{D}_n(f),\quad u\in X\setminus\mathfrak{N},
\quad f\in L^2(X,\mu). \end{equation}

Initially, we will use the above kernels to prove the following
result.

\begin{thm} If $K\in L^2(X\times X,\mu\times \mu)$ and
$n\geq 1$ then $E_n(K)$ is continuous
$\mu\times\mu$-a.e..\end{thm}

\pf It suffices to show that $E_n(K)$ is continuous in the set
$(X\setminus \mathfrak{N})\times (X\setminus \mathfrak{N})$.\ Let
$u_0,v_0\in X\setminus \mathfrak{N}$.\ It is not hard to see that
\begin{equation}\label{3.12}
E_n(K)(u,v)=\int_X\int_X D_n(u,x) K(x,y) D_n(v,y) \,
d\mu(y)d\mu(x),\quad u,v\in X,\end{equation}
and that we can use Lemma 3.1 and Lemma 3.2 to deduce that
\begin{equation}\label{3.13}
\lim_{(u,v)\to (u_0,v_0)}D_n(u,x)K(x,y)D_n(v,y)=
D_n(u_0,x)K(x,y)D_n(v_0,y), \end{equation}
for $x,y \in X$ a.e..\ If $(u,v)\in O_n(u_0)\times O_n(v_0)$, we have
\begin{equation}\label{3.14}
|D_n(u,x)K(x,y)D_n(v,y)|\leq {1\over
\mu(O_n(u_0))\mu(O_n(v_0))}|K(x,y)|,\end{equation}
for $x,y \in X$ a.e..\ So, the continuity at $(u_0,v_0)$ now follows from the dominated
convergence theorem.\eop

Next, we will state and prove a list of technical results that will
lead to the following conclusion: $\mathcal{D}_n\mathcal{K}\mathcal{D}_n$ coincides with the
integral operator generated by $E_n(K)$.

The sequence of partitions $\{\mathcal{P}_n\}$ was constructed in such a way
that each one of them has the following feature: every element of
$\{\mathcal{A}_i\}$ is a subset of at most finitely many $O_n(x)$.\ That been
said, if $n$ and $i$ are fixed, we can write
\begin{equation}\label{3.15}
\mathcal{A}_i\subset \left(\bigcup_{j=1}^{m(n,i)}O_n(x_j)\right)\bigcup
\mathfrak{N}(n,i),\end{equation}
in which $\mu(\mathfrak{N}(n,i))=0$ and $0<\mu(O_n(x_j))<\infty$,
$j=1,2, \ldots, m(n,i)$.\ The set $\mathfrak{N}(n,i)$ is nothing but
the union of all elements of $\mathcal{P}_n$ for which the
intersection with $\mathcal{A}_i$ has measure zero.

In the next results, we will deal with a continuous function $f:
X\to \mathbb{C}$ with compact support $X_f$.\ Since $X_f$ can be covered by
finitely many $\mathcal{A}_i$, after re-ordering if necessary, we can find an
index $l$ so that
\begin{equation}\label{3.16}
X_f \subset
                 \left(\bigcup_{k=1}^{l}\bigcup_{j=1}^{m(n,k)}O_n(x_j)\right)\bigcup
                 \left(\bigcup_{k=1}^{l}\mathfrak{N}(n,k)\right),\end{equation}
with $\mu(\cup_{k=1}^{l}\mathfrak{N}(n,k))=0$.\ In that case, we
will write
\begin{equation}\label{3.17}
Y_f=\bigcup_{k=1}^{l}\bigcup_{j=1}^{m(n,k)}O_n(x_j).\end{equation}

\begin{lem} Let $f : X\to \mathbb{C}$ be a function with compact
support $X_f$ and $K$ an element of $L^1_{loc}(X\times X, \mu\times
\mu)$.\ Then \begin{align*}
\int_{X\times X}\int_{X} D_n(u,x) K(x,y)& D_n(y,z)f(z)\,d\mu(z)\, d(\mu\times \mu)(x,y)\\
&=\int_{X} \int_{X\times X}D_n(u,x) K(x,y) D_n(y,z)f(z)\,d(\mu\times
\mu)(x,y)\,d\mu(z).
\end{align*}
\end{lem}

\pf Pick $M>0$ so that $|f(x)|\leq M$, $x\in X$.\ We have
\begin{align*}
\int_{X\times X} \int_{X} |D_n(u,x) &K(x,y) D_n(y,z)f(z)|\,d\mu(z)\, d(\mu\times \mu)(x,y)\\
&\leq {M\over \mu(O_n(u))} \int_{O_n(u)\times X}{|K(x,y)|\over
\mu(O_n(y))} \int_{X_f} \chi_{O_n(y)}(z)\,d\mu(z)\,d(\mu\times
\mu)(x,y).
\end{align*}
If $y\not\in Y_f$ then $y\not\in X_f$ and, consequently, $O_n(y)\cap
X_f=\emptyset$.\ Thus $\chi_{O_n(y)}=0$ in $X_f$ and we can take the
above integral on $O_n(u)\times Y_f$.\ Now, if $u\in X\setminus
\mathfrak{N}$, the local integrability of $K$ implies that
\begin{align*}
\int_{O_n(u)\times Y_f}{|K(x,y)|\over \mu(O_n(y))} \int_{X_f}
\chi_{O_n(y)}(z)\,d\mu(z)&\,d(\mu\times \mu)(x,y)\\
&= \int_{O_n(u)\times Y_f}{|K(x,y)|\over \mu(O_n(y))}\, \mu(O_n(y)\cap X_f) \,d(\mu\times \mu)(x,y)\\
&\leq \int_{O_n(u)\times Y_f}|K(x,y)| \,d(\mu\times
\mu)(x,y)<\infty.
\end{align*}
Fubini's theorem ([12, p.386]) completes the proof.\eop

\begin{lem} Let $f: X\to \mathbb{C}$ be a continuous function with
compact support $X_f$ and $K$ an element of $L^1_{loc}(X\times X,
\mu\times \mu)$.\ If $u\in X\setminus \mathfrak{N}$ then
\begin{align*}
\int_{O_n(u)\times Y_f} K(x,y) D_n(y,z)\,d(\mu\times\mu)(x,y)&=
\int_{O_n(u)} \int_{Y_f}K(x,y) D_n(y,z)\,d\mu(y)\,d\mu(x)\\
&=\int_{Y_f} \int_{O_n(u)} K(x,y) D_n(y,z)\,d\mu(x)\,d\mu(y),
\end{align*}
with $Y_f$ as defined in $(3.17)$.\end{lem}

\pf If $u\in X\setminus \mathfrak{N}$ and $z\in X$ then
\begin{align*}
\int_{O_n(u)} \int_{Y_f} |K(x,y) D_n(y,z)|\,d\mu(y)\,d\mu(x) &=
\int_{O_n(u)} \int_{Y_f}\left|K(x,y)\,{\chi_{O_n(y)}(z)
\over\mu(O_n(y))}\right|\,d\mu(y)\,d\mu(x)\\
&\leq\int_{O_n(u)} \int_{Y_f}{|K(x,y)| \over
\mu(O_n(y))}\,d\mu(y)\,d\mu(x).
\end{align*}
Introducing the decomposition $(3.17)$ in the last expression above
and recalling the uniqueness property of the $O_n(x)$, we deduce
that
\begin{align*}
\int_{O_n(u)} \int_{Y_f} {|K(x,y)| \over
\mu(O_n(y))}\,d\mu(y)\,d\mu(x) &= \int_{O_n(u)}
\sum_{k=1}^{l}\sum_{j=1}^{m(n,k)}\int_{O_n(y_j)}{|K(x,y)| \over
\mu(O_n(y))}\,d\mu(y)\,d\mu(x)  \\
 &=  \int_{O_n(u)} \sum_{k=1}^{l}\sum_{j=1}^{m(n,k)}{1\over \mu(O_n(y_j))}
 \int_{O_n(y_j)}|K(x,y)| \,d\mu(y)\,d\mu(x)\\
  &\leq \max_{\tiny{\begin{array}{c}1\leq k\leq l \\ 1\leq j\leq m(n,k)\\
  \end{array}}}\left\{{1\over \mu(O_n(y_j))}\right\}\|K\|_{L^1(O_n(u)\times Y_f)}< \infty.
\end{align*}
Once again, Fubini's theorem leads to the concluding statement. \eop

\begin{lem} Let $f: X\to \mathbb{C}$ be a continuous function with
compact support $X_f$ and $K$ an element of $L^1_{loc}(X\times
X,\mu\times \mu)$.\ Then
\begin{align*}
\int_{X}\int_X\int_{X} D_n(u,x) K(x,y) &D_n(y,z)f(z)\,d\mu(z)\, d\mu(y)\,d\mu(x)\\
&= \int_{X\times X}\int_{X}D_n(u,x) K(x,y) D_n(y,z)f(z)\,d\mu(z)
\,d(\mu\times \mu)(y,x).
\end{align*}\end{lem}

\pf If $u\in X\setminus \mathfrak{N}$ and $M>0$ is a bound for $f$
in $X$ then it is easily seen that
\begin{align*}
\int_{X}\int_X | \int_{X}D_n(u,x) K(x,y)& D_n(y,z)f(z)\,d\mu(z)| \, d\mu(y)\,d\mu(x)\\
& \leq {M \over \mu(O_n(u))}\int_{O_n(u)}\int_{Y_f}{|K(x,y)|\over \mu(O_n(y))}\,
\mu(O_n(y)\cap X_f)\,d\mu(x) \,d\mu(y)\\
& \leq {M \over
\mu(O_n(u))}\int_{O_n(u)}\int_{Y_f}|K(x,y)|\,d\mu(x)
\,d\mu(y)<\infty.
\end{align*}
So, the result follows from Fubini's theorem once again. \eop

The proof of the next lemma is analogous and will be
omitted.

\begin{lem} Let $f: X\to \mathbb{C}$ be a continuous function with
compact support and $K$ and element in $L^1_{loc}(X\times
X,\mu\times \mu)$.\ Then
\begin{align*}
\int_{X}\int_X\int_{X} D_n(u,x) K(x,y)& D_n(y,z)f(z)\,d\mu(z)\, d\mu(y)\,d\mu(x)\\
&= \int_{X}\int_{X}\int_{X}D_n(u,x) K(x,y) D_n(y,z)f(z)\,d\mu(x)\,
d\mu(y)\,d\mu(z).
\end{align*}\end{lem}

Recalling that if $x,y \in X$ then $z\in O_n(y)$ if and only if
$y\in O_n(z)$, the following lemma becomes obvious.

\begin{lem}If $n\geq 1$ then $D_n(y,z)=D_n(z,y)$, $y,z\in
X\setminus \mathfrak{N}$.\end{lem}

Below, $\mathcal{E}^n_{K}$ will denote the integral operator generated by
$E_n(K)$.

\begin{thm} If $K\in L^2(X\times X,\mu\times \mu)$ then
$\mathcal{D}_n\mathcal{K}\mathcal{D}_n(f)=\mathcal{E}^n_K(f)$, $f \in L^2(X,\mu)$.\end{thm}

\pf Clearly $L^2(X\times X,\mu\times \mu) \subset L^1_{loc}(X\times
X,\mu\times \mu)$.\ If $K\in L^2(X\times X,\mu\times\mu)$ and
$f\colon X\to \mathbb{C}$ is continuous with compact support then the
previous lemmas imply that
\begin{eqnarray*}
\mathcal{D}_n\mathcal{K}\mathcal{D}_n(f)(u) \hspace*{-3mm} &=& \hspace*{-3mm}  \int_{X}\int_X\int_{X} D_n(u,x) K(x,y) D_n(y,z)f(z)\,d\mu(z)\, d\mu(y)\,d\mu(x)\\
                \hspace*{-3mm}   &=&  \hspace*{-3mm} \int_{O_n(u)}\int_{Y_f}\int_{X_f} \hspace*{-2mm}  D_n(u,x) K(x,y) D_n(y,z)f(z)d\mu(z)d\mu(y)d\mu(x)\\
                \hspace*{-3mm}   &=& \hspace*{-3mm}  \int_{X}\int_{X}\int_{X}D_n(u,x) K(x,y) D_n(z,y)f(z)\,d\mu(x)\, d\mu(y)\,d\mu(z)\\
                \hspace*{-3mm}   &=& \mathcal{E}^n_K(f)(u), \quad u \in X\setminus \mathfrak{N}.
\end{eqnarray*}
Hence, the result in the statement of the theorem follows from the
equality $\mu(X\setminus \mathfrak{N})=0$ and from a basis
approximation theorem from measure theory ([12, p.197]). \eop

The last result of the section refers to the positive definiteness
of $E_n(K)$.

\begin{thm} If $K$ is $L^2(X,\mu)$-positive definite then
so is $E_n(K)$.\end{thm}

\pf If $K$ is $L^2(X,\mu)$-positive definite then both, $K$ and
$E_n(K)$, belong to the space $L^2(X\times X,\mu\times\mu)$.\ On the
other hand, Theorem 2.1 and $(3.11)$ imply that
$\mathcal{D}_n(L^2(X,\mu))\subset L^2(X,\mu)$.\ Thus, an application of
Theorem 3.9 leads to
\begin{equation}\label{3.18}
\langle \mathcal{E}^n_K(f),f \rangle_2
=
\langle \mathcal{K}\mathcal{D}_n(f),\mathcal{D}_n(f) \rangle_2 \geq 0,\quad f\in
L^2(X,\mu).\end{equation}
The proof is complete.\eop

\section{Traceability}

This section contains the main results of the paper.\ They can be
interpreted as generalizations of results obtained in [8] and other
references quoted here.\ The traceability results described here will be obtained via several
known results on trace-class operators and singular values of
operators.\ We will quote some of them and just mention others.\ The construction developed in Section 2 reveals that the diagonal of
$X$ is, up to a set of measure zero, a subset of $(X\setminus
\mathfrak{N}) \times (X\setminus \mathfrak{N})$.\ This remark
justify why some of the integrals appearing below are not
identically zero.\ Given $K \in L^2(X\times X,\mu\times\mu)$, we will
consider $\mathcal{E}^n_K$ acting like an operator on $L^2(X,\mu)$.\ All
other operators mentioned here are to be understood acting in the
same way.

The following lemma is an adapted to our purposes version of Theorem
4.1 in [3].

\begin{lem} If $K$ is a continuous ($\mu\times\mu$-a.e.)
$L^2(X,\mu)$-positive definite kernel and $x \in X \rightarrow
K(x,x)$ is integrable then $\mathcal{K}$ is trace-class and
\begin{equation}\label{4.1}\mbox{tr\,}(\mathcal{K})=\int_X K(x,x)\, d\mu(x).\end{equation}
\end{lem}

\begin{lem} Let $K$ be $L^2(X,\mu)$-positive definite.\ If
\begin{equation}\label{4.2}
\limsup_{n\to\infty}\int_{X}E_n(K)(x,x)\,d\mu(x)<
\infty,\end{equation} then
$\limsup_{n\to\infty}\mbox{tr\,}(\mathcal{E}^n_K)<\infty$.\end{lem}

\pf If $(4.2)$ holds then there exists $n_0\in\mathbb{N}$ such that $x\in
X\to E_n(K)(x,x)$ is integrable for $n\geq n_0$.\ Theorem 3.10
implies that $E_n(K)$ is $L^2(X,\mu)$-positive definite while
Theorem 3.3 shows that $E_n(K)$ is continuous $\mu\times \mu$-a.e..\
Applying Lemma 4.1 we see that
\begin{equation}\label{4.3}
\mbox{tr\,}(\mathcal{E}^n_K)=\int_X E_n(K)(x,x)\,d\mu(x),\quad n\geq n_0.\end{equation}
The result follows. \eop

Next, we recall some facts involving singular values of an
operator.\ If $T$ is a compact operator on a Hilbert space, a {\em
singular value of $T$} is an eigenvalue of $(T^{*}T)^{1/2}$.\ We
shall enumerate the nonzero singular values of $T$ in decreasing
order, taking multiplicities into account: $s_1(T)\geq
s_2(T)\geq \ldots$.\ If the rank $\rho$ of $(T^{*}T)^{1/2}$ is
finite, obviously $s_j(T)=0$, $j \geq \rho+1$.\ If the eigenvalues
of $T$ are ordered like $|\l_1(T)|\geq |\l_2(T)|\geq \ldots$, then a
classical result from operator theory states that
$s_j(T)=|\l_j(T)|$, $j=1,2,\ldots$, as long as $T$ is either
hermitian or normal.\ If $S$ is another compact operator of same
type as $T$, and assuming the same ordering on the singular values
of $S$, the following inequality holds: $|s_n(T)-s_n(S)|\leq
\|T-S\|$, $n=1,2\ldots$.\ All of these results can be found with
proofs in $[10,11]$.

In Theorem 4.4, a complement of Lemma 4.2, we also use the following
nontrivial result on convergence of operators ([8]).

\begin{lem} Let $\{T_n\}$ be a countable set of bounded linear
operators on a Hilbert space $\mathcal{H}$ such that
$\lim_{n\to\infty}\|T_n(f)-f\|_{\mathcal{H}}=0$, $f\in\mathcal{H}$.\ If every $T_n$ is
self-adjoint and $T$ is a bounded compact operator on $\mathcal{H}$ then
$\lim_{n\to\infty}\|T_nTT_n-T\|=0$.\end{lem}

\begin{thm} Let $K$ be $L^2(X,\mu)$-positive definite.\ If
\begin{equation}\label{4.4}
\limsup_{n\to\infty}\int_{X}E_n(K)(x,x)\,d\mu(x)<\infty,
\end{equation} then $\mathcal{K}$ is trace-class.\end{thm}

\pf Since $\{s_j(\mathcal{E}^n_K)\}\subset(0,\infty)$, it is quite clear that
\begin{equation}\label{4.5}
\sum_{j=1}^k s_j(\mathcal{E}^n_K)\leq\mbox{tr\,}(\mathcal{E}^n_K),\quad
k=1,2,\dots.\end{equation} Theorem 3.9 and the inequality mentioned
before Lemma 4.3 imply that \begin{equation}\label{4.6}
|s_j(\mathcal{E}^n_K)-s_j(\mathcal{K})|\leq \|\mathcal{E}^n_K-\mathcal{K}\|=\|\mathcal{D}_n\mathcal{K}\mathcal{D}_n-\mathcal{K}\|,\quad
j=1,2,\dots.\end{equation} Since each $\mathcal{D}_n$ is self-adjoint, $\mathcal{K}$ is
compact and
\begin{equation}\label{4.7}
\lim_{n\to\infty}\|\mathcal{D}_n(f)-f\|_2=0,\quad f\in L^2(X,\mu),\end{equation}
we are authorized to apply Lemma 4.3 to conclude, from $(4.6)$, that
\begin{equation}\label{4.8}
\lim_{n\to\infty}s_j(\mathcal{E}^n_K)=s_j(\mathcal{K}), \quad j=1,2,\dots.
\end{equation} It is now clear that
\begin{equation}\label{4.9}
\sum_{j=1}^k s_j(\mathcal{K})=\limsup_{n\to\infty}\sum_{j=1}^k
s_j(\mathcal{E}^n_K)\leq \limsup_{n\to\infty}\mbox{tr\,}(\mathcal{E}^n_K), \quad
k=1,2,\dots,\end{equation}
and that concludes the proof. \eop

In order to deal with the converse of the previous result, we will
need the following result ([10, p.51]): if $S_1$, $S_2$ and $T$ are
bounded linear operators on a Hilbert space and $T$ is compact then
so is the composition $S_1TS_2$ and $s_j(S_1TS_2)\leq
\|S_1\|s_j(T)\|S_2\|$, $j=1,2,\dots$.

\begin{lem} Let $p \in [1,\infty)$ and $K\in L^p(X\times
X,\mu\times\mu)$.\ If $x\in X\to K(x,x)$ is integrable and
$\mu(X)<\infty$ then there is a positive integer $n_0$ for which
$x\in X\to E_n(K)(x,x)$ is integrable when $n\geq
n_0$.\end{lem}

\pf Since \begin{equation}\label{4.10} |E_n(K)(u,u)|\leq
|E_n(K)(u,u)-K(u,u)|+|K(u,u)|,\quad u\in X\setminus \mathfrak{N},
\end{equation}
we can use Doob's convergence theorem to select a positive integer
$n_0$ so that
\begin{equation}\label{4.11}
|E_n(K)(u,u)|\leq 1+|K(u,u)|,\quad  u\in X\setminus \mathfrak{N},
\quad n\geq n_0.\end{equation}
Our assumptions on $X$ and $x\in X\to K(x,x)$ imply the result.
\eop

\begin{thm} Let $K$ be $L^2(X,\mu)$-positive definite.\ If
$x\in X\to K(x,x)$ is integrable and $\mu(X)<\infty$ then there is
$n_0\in\mathbb{N}$ so that $\mathcal{E}_K^n\in \mathcal{B}_1(L^2(X))$ and
\begin{equation}\label{4.12}\mbox{tr\,}(\mathcal{E}_K^n)=\int_X E_n(K)(x,x)\,d\mu(x),\quad n\geq
n_0.\end{equation}\end{thm}

\pf The previous lemma reveals that $x\in X\to E_n(K)(x,x)$ is
integrable for $n$ large.\ As so, the result follows from Theorem
3.3 and Lemma 4.1. \eop

\begin{thm} Let $K\in L^2(X\times X,\mu\times\mu)$.\ If
$\mathcal{K}$ is trace-class then so is every $\mathcal{E}^n_K$.\ The number $\mbox{tr\,}(\mathcal{K})$
is an upper bound for the sequence $\{\mbox{tr\,}(\mathcal{E}^n_K)\}$.\end{thm}

\pf Assume $\mathcal{K}$ is trace-class.\ Since each $\mathcal{D}_n$ is bounded,
Theorem 3.10 and the comments preceding Lemma 4.5 imply that
\begin{equation}\label{4.13}
s_j(\mathcal{E}^n_K)=s_j(\mathcal{D}_n\mathcal{K}\mathcal{D}_n)\leq \|\mathcal{D}_n\| s_j(\mathcal{K})\|\mathcal{D}_n\|,\quad
n=1,2,\dots.\end{equation}
Hence,
\begin{equation}\label{4.14}
\sum_{j=1}^{\infty}s_j(\mathcal{E}^n_K)\leq
\|\mathcal{D}_n\|^2\sum_{j=1}^{\infty}s_j(\mathcal{K}),\end{equation}
and the two assertions of the lemma follow.\eop

The following result is very close to a converse of Theorem
4.4.

\begin{thm} Let $K\in L^2(X\times X,\mu\times\mu)$.\ If
$\mathcal{K}$ is trace-class then $$\lim_{n\to\infty}
\mbox{tr\,}(\mathcal{E}^n_K)=\mbox{tr\,}(\mathcal{K}).$$\end{thm}

\pf A basic inequality for the trace ([10, p.54]) implies that
\begin{equation}\label{4.15}
|\mbox{tr\,}(\mathcal{E}^n_K)-\mbox{tr\,}(\mathcal{K})|\leq \sum_{j=1}^{\infty}s_j(\mathcal{E}^n_K-\mathcal{K}),\quad
n=1,2,\ldots,\end{equation} as long as $\mathcal{K}$ is trace-class.\ On the other hand, since (see [10, p.89])
\begin{equation}\label{4.16}\lim_{n\to \infty}\sum_{j=1}^{\infty}s_j(\mathcal{D}_n \mathcal{K} \mathcal{D}_n-\mathcal{K})=0,\end{equation}
Theorem 3.9 completes the proof.\eop

Next, we move to a proof of the converse of Theorem 4.3 in the case
when $\mu(X)<\infty$.

\begin{thm} Let $K$ be $L^2(X,\mu)$- positive definite.\ If
$\mathcal{K}$ is trace-class and $\mu(X)<\infty$ then \begin{equation}\label{4.17}\lim_{n\to\infty}
\int_X E_n(K)(x,x)\,d\mu(x)<\infty.\end{equation}\end{thm}

\pf Assume $\mathcal{K}$ is trace-class.\ Since the function $x\in X\to
K(x,x)$ is integrable already, if $\mu(X)<\infty$, we can use
Theorem 4.6 to find a positive integer $n_0$ such that
\begin{equation}\label{4.18}tr(\mathcal{E}^n_K)=\int_XE_n(K)(x,x)\,d\mu(x),\quad n\geq n_0.\end{equation} An
application of Theorem 4.8 finishes the proof. \eop

At this point, it is very important to remind the reader that the
results we have obtained includes the case in which $X$ is either a
sphere or a torus.

Next, we intend to consider cases in which $X$ has no finite
measure.\ In order to handle that, we use the cover $\{\mathcal{A}_m\}$ of
$X$ constructed before to define a sequence of subsets of $X$ that
increases to $X$.\ Precisely, defining $X_j=\cup_{m=1}^{j}\mathcal{A}_m$,
$j\geq 1$, we immediately have the following two properties:
$X=\cup_{j=1}^{\infty}X_j$ and if $x\in X$ then there exists $j_0
\geq 0$ such that $x\in X_{j}$, $j\geq j_0$.\ Using the sequence
just defined, we now take linear operators $P_j:
 L^2(X,\mu)\to L^2(X,\mu)$ defined by the formula $P_j(f)=f\chi_{X_j}$, $f\in
L^2(X,\mu)$.\ They are self-adjoint and the uniform boundedness
principle shows that the sequence $\{P_j\}$ is bounded in the space
of bounded linear operators on $L^2(X,\mu)$.\ Also, the dominated
convergence theorem implies that $\{P_j\}$ converges pointwise to
the identity operator on $L^2(X,\mu)$.\ The following technical
lemma contains a critical information on the sequence
$\{P_j\}$.

\begin{lem} If $T: L^2(X,\mu) \rightarrow L^2(X,\mu)$ is
trace-class then each $P_jTP_j$ is so and the limit formula
$\lim_{j\to\infty}\mbox{tr\,}(P_jTP_j)=\mbox{tr\,}(T)$ holds.\end{lem}

\pf The first assertion is a consequence of the remark preceding
Lemma 4.5.\ As for the other, it follows from Theorem 11.3 in [10]
.\eop

The converse of Theorem 4.4 reads as follows.

\begin{thm} Let $K$ be $L^2(X,\mu)$-positive definite.\ If
$\mathcal{K}$ is trace-class then the limit
\begin{equation}\label{4.19}
\lim_{n\to\infty}\int_XE_n(K)(x,x)\,d\mu(x)\end{equation}
exists and is finite.\end{thm}

\pf The proof requires the double-indexed operator $\mathcal{Q}^n_j$ given by
the formula
\begin{equation}\label{4.20}
\mathcal{Q}^n_j(f)(x)=\int_{X_j}E_n(K)(x,y)f(y)\,d\mu(y),\quad x\in X_j,\quad
f\in L^2(X_j,\mu).\end{equation}
If $f\in L^2(X_j)$, let us write $\tilde{f}$ to denote a function on
$X$ that coincides with $f$ on $X_j$ and is zero in $X\setminus
X_j$.\ It is now clear that
\begin{eqnarray*}
\int_{X_j}\mathcal{Q}^n_j(f)(x)\overline{f(x)}\,d\mu(x) &=& \int_X\int_X
E_n(K)(x,y)
\tilde{f}(y)\overline{\tilde{f}(x)}\,d\mu(y)\,d\mu(x)\\
                  &=&  \int_X\int_X E_n(K)(x,y)\tilde{f}(y)\,d\mu(y)
                  \overline{\tilde{f}(x)}\,d\mu(x)\\
                  &=&  \int_X \mathcal{E}^n_K(\tilde{f})(x)\overline{\tilde{f}(x)}\,d\mu(x),
                  \quad f\in L^2(X_j,\mu).
\end{eqnarray*}
Since $K$ is $L^2(X,\mu)$-positive definite, Theorem 3.10 guarantees
that $\mathcal{Q}^n_j$ is $L^2(X_j,\mu)$-positive definite.\ Also, the fact
that $\mathcal{K}$ is trace-class implies that $x\in X\to K(x,x)$ is
integrable.\ Hence, due to Lemma 4.5, there exists $n_0\geq 0$ such
that $x\in X_j\to E_n(K)(x,x)$ whenever $n\geq n_0$.\ Recalling
Theorem 3.3 and applying Lemma 4.1, we deduce that $\mathcal{Q}^n_j$ is
trace-class and
\begin{equation}\label{4.21}
\mbox{tr\,}(\mathcal{Q}^n_j)=\int_{X_j}E_n(K)(x,x)\,d\mu(x),\end{equation}
as long as $n\geq n_0$.\ Let us keep the previous condition on $n$
in force.\ If $V_j$ is the closed subspace $L^2(X,\mu)$ encompassing
the functions on $X$ which are zero in $X\setminus X_j$ and $\mathcal{R}^n_j
: V_j\to V_j$ is the operator given by
\begin{equation}\label{4.22}
\mathcal{R}^n_j(f)(x)=\chi_{X_j}(x)\int_{X}E_n(K)(x,y)\chi_{X_j}(y)f(y)\,d\mu(y),
\end{equation} with $x\in X, f\in V_j$,
then $\mathcal{R}^n_j$ and $\mathcal{Q}^n_j$ possess the same eigenvalues.\ Having in
mind the previous lemma,
\begin{equation}\label{4.23}
(\mathcal{P}_j\mathcal{E}^n_K\mathcal{P}_j)(f)(x)=\int_X E_n(K)(x,y)\chi_{X_j\times
X_j}(x,y)f(y)\,d\mu(y),\end{equation} for $x\in X, f\in L^2(X)$,
and we can conclude now that $ \mathcal{R}^n_j$ and $\mathcal{P}_j\mathcal{E}^n_K\mathcal{P}_j$ have the
same eigenvalues.\ Therefore, \begin{equation}\label{4.24}
\mbox{tr\,}(\mathcal{P}_j\mathcal{E}^n_K\mathcal{P}_j)=\mbox{tr\,}(\mathcal{R}^n_j)=\mbox{tr\,}(\mathcal{Q}^n_j)=\int_{X_j}E_n(K)(x,x)\,d\mu(x).\end{equation}
The monotone convergence theorem leads to
\begin{equation}\label{4.25}
\mbox{tr\,}(\mathcal{E}^n_K)=\int_X E_n(K)(x,x)\,d\mu(x).\end{equation}
Finally, (4.24) and the observation made before Theorem 4.3 lead to
the assertion of the theorem.\eop

\vspace*{6mm} \footnotesize{\noindent Departamento de Matem\'{a}tica,\\
\noindent ICMC-USP - S\~{a}o Carlos, Caixa Postal 668,\\
13560-970 S\~{a}o Carlos SP, Brasil. \\
menegatt@icmc.usp.br\\
apperon@icmc.usp.br}

\vspace*{6mm} \footnotesize{\noindent FAMAT-UFU \\Caixa Postal 593\\
38400-902 Uberl\^{a}ndia-MG, Brasil\\
mariocastro@famat.ufu.br}

\end{document}